\documentclass[a4paper,11pt]{amsart}
\usepackage[left=2.5cm, right=2.5cm,top=2.5cm,bottom=3cm]{geometry}
\usepackage{amssymb}
\usepackage{graphicx}
\usepackage{setspace,algorithmic}

\usepackage{subcaption}
\usepackage[utf8]{inputenc}
\usepackage[english]{babel}
\usepackage{amsfonts}
\usepackage{dsfont}

\usepackage{hyperref}
\usepackage{enumitem}
\usepackage{tikz}
\usepackage{chemarr}

\usepackage{graphicx}
\usepackage{hhline}

\usepackage{tikz-cd}
\usepackage{algorithm}
\usepackage{algorithmic}

\newtheorem{thm}{Theorem}[section]
\newtheorem{cor}[thm]{Corollary}
\newtheorem{lem}[thm]{Lemma}
\newtheorem{prop}[thm]{Proposition}

\theoremstyle{dfn}

\theoremstyle{exa}
\newtheorem{exa}[thm]{Example}
\theoremstyle{rem}
\newtheorem{rem}[thm]{Remark}
\numberwithin{equation}{section}

\newcommand{\R}{\mathbb{R}}
\newcommand{\Rs}{\R\setminus\{0\}}
\newcommand{\Rn}{\mathbb{R}^n}
\newcommand{\eps}{\varepsilon}

\newcommand{\cA}{\mathcal{A}}
\newcommand{\hA}{\widehat{\cA}}
\newcommand{\cH}{\mathcal{H}}

\DeclareMathOperator*{\Log}{Log}

\DeclareMathOperator*{\diag}{Diag}

\DeclareMathOperator{\Conv}{Conv}

\DeclareMathOperator{\Sing}{Sing}

\DeclareMathOperator{\sign}{sign}
\DeclareMathOperator{\Hess}{Hess}

\newcommand{\floor}[1]{\left\lfloor#1\right\rfloor}

\renewcommand{\qed}{$\blacksquare$}
\newcommand{\dia}{$\diamond$}

\begin{document}

\title{Optimal Bounds for the Number of Pieces of Near-Circuit
Hypersurfaces}

\author{Weixun Deng}
\email{deng15521037237@tamu.edu }

\author{J.\ Maurice Rojas}
\email{rojas@tamu.edu}

\author{Cordelia Russell}
\email{cordelia03@gmail.com}

\renewcommand{\shortauthors}{Weixun Deng, J.\ Maurice Rojas and Cordelia Russell}

\begin{abstract}
Suppose $f$ is a polynomial in $n$ variables with real coefficients, exactly
$n+k$ monomial terms, and Newton polytope of positive volume. Estimating the
number of connected components of the positive zero set of $f$ is a fundamental
problem in real algebraic geometry, with applications in computational
complexity and topology. We prove that the number of connected components is
at most $3$ when $k\!=\!3$, settling an open question from Fewnomial Theory.
Our results also extend to exponential sums with real exponents. A key
contribution here is a deeper analysis of the underlying $\cA$-discriminant
curves, which should be of use for other quantitative geometric
problems.
\end{abstract}

\keywords{Exponential Sum, Morse Theory, Connected Component, Discriminant}

\maketitle

\section{Introduction}
Estimating the number of connected components (a.k.a. {\em pieces}) of the
real zero set of a polynomial is a fundamental problem with numerous
applications. For univariate polynomials with exactly $t$ monomial terms,
Descartes' Rule of Signs reveals that the maximal number of isolated positive
roots is at most $t-1$. In the 1980s, Khovanskii
\cite{Khovanskii} extended this to higher dimensions by developing {\em
Fewnomial Theory}. One of his bounds gave, for any $n$-variate polynomial $f$
with real coefficients and exactly $n+k$ monomial terms, an upper bound of
$2^{O((n+k-1)^2)}$ for the number of pieces of the zero
set, $Z_+(f)$, of $f$ in the positive orthant $\Rn_+$.\footnote{In our use of
$O$- and $\Omega$- notation, all our constants are effective and absolute,
i.e., they can be made explict, and there is {\em no} dependence on any further
parameters.}

Although Khovanskii's fewnomial bounds were eventually found to be far from
optimal, significant improvements took decades to find: After improvements by
Li, Rojas, Wang \cite{LiRojasWang::FirstBound}, and Perrucci
\cite{Perrucci::PerrucciBound}, Bihan and Sottile
sharpened the bounds to
$2^{O(k^2+n+k\log n)}$
\cite{BihanRojasSottile::BRSBound,BihanSottile::BSBound}.

Refining an approach from \cite{forsgard2017new}, Bihan, Humbert,
and Tavenas \cite{BihanBound} then applied \(\mathcal{A}\)-discriminants
\cite{gelfand1994discriminants} to prove
an even sharper upper bound of the form $2^{O(k^2+k\log n)}$.
In particular, when $k\!=\!3$ and the exponent vectors do {\em not} lie
in an affine hyperplane, the underlying $f$ is called an {\em honest}
(real) $n$-variate $(n+3)$-nomial, or a {\em near-circuit polynomial}, and
their best upper bound for this case was
$\floor{\frac{n-1}{2}} + 3$. We also call the underlying set of
exponent vectors a {\em near-circuit}.

\begin{thm}
\label{Thm::Main}
Let $f$ be an honest real $n$-variate $(n+3)$-nomial. Then the number of
pieces of $Z_+(f)$ is at most 3.
\end{thm}

\noindent
The following example shows that our bound is in fact optimal.
\begin{exa}
\label{Ex::Optimal}
\cite[Ex.\ 1.8]{prt}
Consider\\
\mbox{}\hfill $f(x_1,x_2)\!:=\!1-x_1-x_2+\frac{6}{5}x_1x^4_2
+\frac{6}{5}x^4_1x_2$.\hfill\mbox{}\\
The coordinate-wise natural log map of $Z_+(f)$ is drawn below:\\
\mbox{}\hfill \includegraphics[scale=0.35]{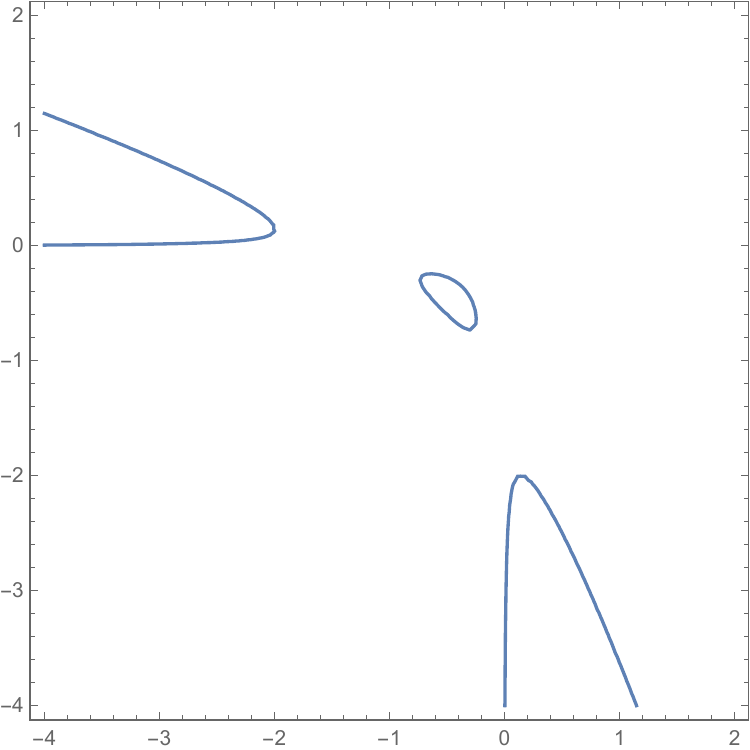}
\hfill\mbox{}\\
A standard computation of critical points of the projection mapping
$Z_+(f)$ to $\R^{n-1}_+\times \{0\}$ easily yields a proof that we indeed
have $3$ pieces, exactly one of which is bounded. \dia
\end{exa}
\begin{rem}
The correct maximal numbers of pieces in the honest $n$-variate
$(n+k)$-nomial case, for $k\!\in\{1,2\}$, are respectively $1$ and $2$. It is
worth noting that it took about $18$ years from the publication of Khovanskii's
book on fewnomials \cite{Khovanskii} until the 2011 habilitation thesis of
Bihan for the second upper bound to be settled. The maximal number of
pieces for honest $n$-variate $(n+4)$-nomials remains unknown. \dia
\end{rem}

To achieve our bound we closely analyze the structure of
{\em \(\mathcal{A}\)-discriminants} for near-circuit polynomials
(see Subsections \ref{A-disc for n+3} and \ref{Chambers Structure}).
The key is that once a near-circuit $\cA$ is fixed, the resulting
$(n+3)$-dimensional family of near-circuit polynomials can be parametrized
via the pieces of the complement of a (possibly singular) algebraic curve
$\nabla\!\subset\!\R^2_+$. The number of cusps on this (reduced) discriminant
curve determines, in part, the number of possible isotopy types for
$Z_+(f)$, building on earlier work of
work of Rojas and Rusek \cite{RojasRusek_Adiscriminant} (see also Lemmas \ref{Thm::RojasRusek_Adiscriminant} and \ref{Lem::diff of adjacent chambers}).
The more challenging cases arises when the number of cusps is two or more,
particularly because we need to determine how the isotopy type of $Z_+(f)$
changes as we vary the coefficients of $f$ and traverse $\nabla$.

Our investigation of how cusps on $\cA$-discriminant curves affect the possible
isotopy types of $Z_+(f)$ is a tool that we hope will be used more broadly in
real algebraic geometry.

\section{Preliminaries}
For any two vectors $v\!=\!(v_1,\ldots,v_n),w\!=\!(w_1,\ldots,w_n)\!\in\!\Rn$
we set $v\cdot w\!:=\!v_1w_1+\cdots+v_nw_n$.
The transpose of a matrix $M$ will be denoted by $M^\top$. For
any function $h: \mathbb{C}^n\to \mathbb{R}$, we set $\Rs\!:=\!\R\setminus
\{0\}$ and let $Z_{\mathbb{R}}(h)$
denote the zero set of $h$ in $\mathbb{R}^n$.

\subsection{Deforming Zero Sets of Exponential Sums}
Suppose $\mathcal{A} = \{\alpha_1, \dots , \alpha_{n+k} \} \subseteq
\mathbb{R}^{n}$ has cardinality $n+k$. Our results in fact hold in the
more general context of real exponents. So let
$f_c(x):= \sum_{i=1}^{n+k} c_i e^{\alpha_i\cdot x}$,
where $c \in (\Rs)^{n+k}$. Since $\log$ defines
a homemorphism from $\R$ to $\R_+$, we will thus henceforth focus on
bounding the number of pieces of $Z_+(f_c)$. (In fact, all the fewnomial
bounds we have mentioned earlier hold at this level of generality.)

We call
$\varepsilon = \sign( c) \in \{ \pm 1 \}^{n+k} $ a \emph{sign distribution},
and $(\mathcal{A},\varepsilon)$ a \emph{signed support}. We also say that
$f_c$ is \emph{honestly $n$-variate} iff the dimension of
$\Conv(\mathcal{A})$ is $n$. We write $\R^{n+k}_\varepsilon = \{ c \in
\R^{n+k} \; | \; \sign(c) = \varepsilon \}$
for the appropriate sub-orthant of $\mathbb{R}^{n+k}$.

Two subsets $Z_0, Z_1 \subseteq \mathbb{R}^n$ are \emph{isotopic} (ambiently
in $\mathbb{R}^n$) iff there is a continuous map $H\colon [0,1] \times
\mathbb{R}^n \to \mathbb{R}^n$ such that
(1) $H(t, \cdot)$ is a homeomorphism for all $t \in [0,1]$,
(2) $H(0,\cdot)$ is the identity on  $\mathbb{R}^n$,
and (3) $H(1,Z_0) = Z_1$.
Isotopy is an equivalence relation on subsets of
$\mathbb{R}^n$ \cite[Ch.\ 10.1]{Book::Isotopy}. So we can speak of
isotopy {\em type}.

\subsection{Signed $\mathcal{A}$-discriminant Contours}
\label{Sec::SignedDiscr}
We recall the notion of \emph{$\mathcal{A}$-discriminant} from
\cite{gelfand1994discriminants, BihanBound}, but
extended to real exponents as in \cite{RojasRusek_Adiscriminant}. A point $x \in \mathbb{R}^n$ is a \emph{singular zero} of $f_c$ ifif
$f_c(x) = \frac{\partial f_c(x) }{\partial x_1} = \dots =
\frac{\partial f_c(x) }{\partial x_n} = 0$.
We denote the set of singular zeros of $f_c$ by $\Sing(f_c)$. For a
fixed signed support $(\mathcal{A},\varepsilon)$, we define the
\emph{signed $A$-discriminant variety} as
\begin{align*}
    \nabla_{\mathcal{A},\varepsilon} := \big\{ c \in \mathbb{R}^{n+k}_\varepsilon \mid \Sing(f_c) \neq \emptyset \big\}.
\end{align*}

We recall a natural invariance property of the signed
$\mathcal{A}$-discriminant.
\begin{prop}
\cite{DengRojasTelek}
\label{Lemma::Transform}
Let $f_c$ be an exponential sum with support
$\mathcal{A} = \{ \alpha_1, \dots , \alpha_{n+k} \} \subseteq \mathbb{R}^n$.
For an invertible matrix $M \in \mathbb{R}^{n \times n}$ and
$v \in \mathbb{R}^n$ consider the exponential sum
$g_c(x) = \sum_{i = 1}^{n+k} c_i e^{ (M \alpha_i  + v)  \cdot x}$.
Then we have:
\begin{itemize}
    \item[(i)]  If $\det(M) > 0$, then the hypersurfaces $Z_{\mathbb{R}}(f_c)$ and $Z(g_c)$ are isotopic.
    \item[(ii)] $\Sing(f_c) = M^{\top}\Sing(g_c)$.
    \item[(iii)] For all $x \in  \Sing(g_c)$ the Hessian matrices $\Hess_{f_c}(M^{\top}x)$ and $\Hess_{g_c}(x) $ have the same number of positive, negative and zero eigenvalues.
\end{itemize}
\end{prop}

\begin{rem}
Via Proposition \ref{Lemma::Transform} one can transform any full-dimensional
support $\mathcal{A} \subseteq \mathbb{R}^n$ to a support containing the
standard basis vectors $e_1, \dots, e_n \in \mathbb{R}^n$ and the zero vector
without changing the isotopy types of the corresponding hypersurfaces.
\end{rem}

For any support $\mathcal{A} = \{ \alpha_1, \dots , \alpha_{n+k} \} \subseteq
\mathbb{R}^n$, let $\hA\!\in\!\R^{(n+1)\times (n+k)}$
denote the matrix with top row $[1,1,\ldots,1]$ and bottom $n$ rows formed from $\mathcal{A}$ by considering each $\alpha_i$ a column vector.
If $x \in \mathbb{R}^n$ is a singular zero of $f_c$, then
$[c_1e^{\alpha_1 \cdot x},\ldots,c_{n+k}e^{\alpha_{n+k}\cdot x}]
\in \ker(\hA)$. Choose a
basis of $\ker(\hat{\mathcal{A}})$ and write these vectors as columns of a
matrix $B \in \mathbb{R}^{(n+k) \times (k-1)}$. Such a choice of $B$ is called
a \emph{Gale dual matrix} of $\mathcal{A}$.

We also call $\mathcal{A}$ \emph{degenerate} iff $\hA$ has rank $\leq\!n$.
For example, when $n=2$ and $k=3$,
\includegraphics[width=0.035\textwidth]{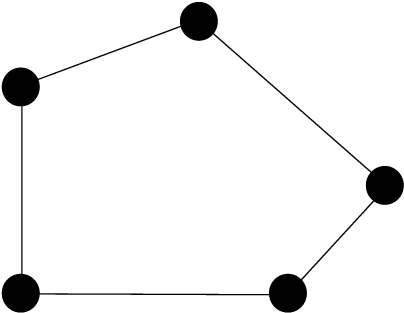} is non-degenerate,
whereas \includegraphics[width=0.035\textwidth]{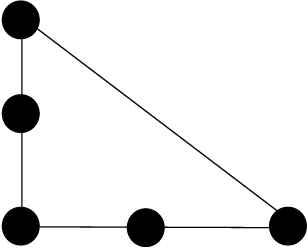} and
\includegraphics[width=0.035\textwidth]{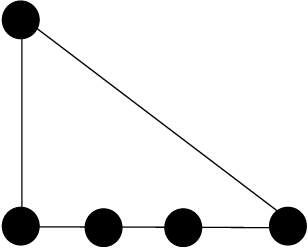} are degenerate.

To understand how the number of pieces of $Z_\R(f_c)$ changes as we vary
$c$, let us recall variants of \cite[Theorem 3.8.]{RojasRusek_Adiscriminant}
and \cite[Theorem 3.14]{forsgard2017new} (and elaborate further in
Section \ref{Morse}):
\begin{lem}
\label{Thm::RojasRusek_Adiscriminant}
Let $(\mathcal{A},\varepsilon)$ be a full-dimensional signed support with Gale dual matrix $B$ and let $c,c' \in \mathbb{R}^{n+k}_\varepsilon$. If
$B^{\top}\Log\lvert c \rvert $ and $B^{\top}\Log\lvert c' \rvert $ are in the
same piece of
\[ \mathbb{R}^{k-1} \setminus \big( \bigcup_{F \subseteq \Conv( \mathcal{A} )  \text{ a face}}  B^{\top}\Log\lvert \nabla_{\mathcal{A}_F,\varepsilon_F} \rvert \big), \]
then the zero sets $Z_{\mathbb{R}}(f_c)$ and $Z_{\mathbb{R}}(f_{c'})$ are  ambiently isotopic in $\mathbb{R}^{n}$.
\end{lem}
\begin{lem}
 \label{Lem::diff of adjacent chambers}
    Under the same assumptions of Lemma \ref{Thm::RojasRusek_Adiscriminant}, if $B^{\top}\Log\lvert c \rvert $ and $B^{\top}\Log\lvert c' \rvert $ lie in
adjacent pieces of
\[ \mathbb{R}^{k-1} \setminus \big( \bigcup_{F \subseteq \Conv( \mathcal{A} )
\text{ a face}}  B^{\top}\Log\lvert \nabla_{\mathcal{A}_F,\varepsilon_F}
\rvert \big), \]
then the difference in the number of pieces of $Z_{\mathbb{R}}(f_c)$ and $Z_{\mathbb{R}}(f_{c'})$ is at most 1.
\end{lem}

Pieces of discriminant complements as above are called {\em chambers}.
Following \cite{RojasRusek_Adiscriminant}, we can reduce the dimension of
our chambers by quotienting out some homogeneities, without losing essential
information: We define the
\emph{(signed reduced) $\mathcal{A}$-discriminant contour}
$\Gamma_\varepsilon(\mathcal{A},B)$
\cite[Definition 2.5]{RojasRusek_Adiscriminant} to be
$\Gamma_\varepsilon(\mathcal{A},B) := B^{\top}\Log\lvert \nabla_{\mathcal{A},\varepsilon} \rvert$,
where $\Log$ is the coordinate-wise natural logarithm map and $\lvert \cdot \rvert$ denotes the coordinate-wise absolute value map.

When $\mathcal{A}$ is non-degenerate we have
\[
\bigcup_{F \subseteq \Conv( \mathcal{A} )  \text{ a face}}  B^{\top}\Log\lvert \nabla_{\mathcal{A}_F,\varepsilon_F} \rvert=\Gamma_\varepsilon(\mathcal{A},B).
\]
In \cite{RojasRusek_Adiscriminant}, bounded (resp. unbounded) pieces of $\mathbb{R}^{k-1} \setminus \Gamma_\varepsilon(A,B)$ were called
{\em inner (resp.\ outer) chambers}. The signed
reduced $A$-discriminant $\Gamma_\varepsilon(\mathcal{A},B)$ then admits a
parametrization --- known as the {\em Horn-Kapranov Uniformization} ---
with highly restricted normal vectors:
\begin{prop}
Let $(\mathcal{A},\varepsilon)$ be a non-degenerate signed support with Gale
dual matrix $B$ and set\\
\mbox{}\hfill $\cH_B\!:=\!\{\lambda\!\in\!\R^{k-1}\; | \;
\text{some coordinate of } B\lambda \text{ is } 0\}$\hfill\mbox{}\\
and $C_{B,\eps}
\!:=\!\{\lambda\!\in\!\R^{k-1}\; | \; B\lambda\!\in\!\R^{n+k}_\eps\}$. Also set
$\xi_{B,\varepsilon}(\lambda)\!:=\!B^{\top} \Log \lvert B\lambda \rvert$
for any $\lambda\!\in\!C_{B,eps}$. Then
$\Gamma_\eps(\cA,B)\!=\!\xi_{B,\eps}(C_{B,\eps})$.
\end{prop}
\begin{lem}
\cite[Theorem 2.1]{Kapranov1991ACO}
\label{Lem::Normal Vectors} Any $\lambda\!\in\!\R^{k-1}\setminus\cH_B$
is normal to the point
$\xi_{B,\varepsilon}(\lambda)$, i.e.,
$\sum_{i=1}^{n+k}\lambda_i\cdot\frac{\partial (\xi_{B,\varepsilon})_i}
{\lambda_j}=0$ for all $j$.
\end{lem}

Modifying a Gale dual matrix using elementary column operations gives another choice of Gale dual matrix. Thus, one can assume without any restriction that the last row of the Gale dual matrix $B$ has the form $B_{n+k} = (0,\dots,0,-1)$.
In which case, since $\xi_{B,\varepsilon}$ is homogeneous, one can replace
$\mathbb{R}^{k-1}\setminus\cH_B$ by the $(k-2)-$dimensional quasi-affine
subspace $\{ \lambda \in \mathbb{R}^{k-1}\setminus\cH_B \mid \lambda_{k-1}= 1
\}$. In Section \ref{Morse}, we will prefer this latter choice and work with
a refined version of $\xi_{B,\eps}$ defined as
$\bar{\xi}_{B,\eps}(\mu)\!:=\!B^\top\Log|B[\mu, 1]^\top|$, for any
$\mu\!\in\!\R^{k-2}$ with $B[\mu,1]^\top\!\in\!\R^{n+k}_\eps$.

\subsection{Near-Circuit Exponential Sums}
\label{A-disc for n+3}
Now let us consider the case when $k=3$. Then $B$ will have $2$ columns
and $C_{B,\eps}$ will be an open sub-interval of $\R$ with endpoints
$\displaystyle{-\frac{b_{i,2}}{b_{i,1}}}$ and $\displaystyle{-\frac{b_{j,2}}{b_{j,1}}}$ for some $1\leq i,j\leq n+3$. By changing the order of $\{\alpha_1,\ldots,\alpha_{n+3}\}$, one can assume $i=n+2$ and $j=n+3$. By Proposition \ref{Lemma::Transform} and translation of $x$, we can then reduce $f_c$ to the following
form
\begin{align}
\label{Eq::ReducedExpSum}
\varepsilon_0+\varepsilon_{1}e^{x_1}+\cdots+\varepsilon_{n}e^{x_n}+\varepsilon_{n+1}e^{\beta\cdot x-c_1}+\varepsilon_{n+2}e^{\gamma\cdot x-c_2}
\end{align}
By column operations on $B$, we can then reduce the last two rows of $B$ to
$\begin{bmatrix}
    -1 & 0\\
    0 & -1
\end{bmatrix}$. Finally, since $\hat{\mathcal{A}}B=0$, we then obtain that
the Gale dual matrix $B$ has the following form:
\begin{align}
\label{Eq::ReducedGualDualMatrix}
B=\begin{bmatrix}
1-\sum_{i=1}^{n}\beta_{i} & \beta_1 & \cdots & \beta_{n} & -1 & 0 \\
1-\sum_{i=1}^{n}\gamma_{i} & \gamma_1 & \cdots & \gamma_{n} & 0 & -1
\end{bmatrix}^\top
\end{align}

By the Horn-Kapranov Uniformization $(c_1,c_2)\!=\!\xi_{B,\eps}([\lambda_1,
\lambda_2]$ provided
\begin{align*}
  (c_1,c_2) =(\sum_{i=1}^{n+3} b_{i,1}\log|b_{i,1}\lambda_1+b_{i,2}\lambda_2|, \sum_{i=1}^{n+3} b_{i,2}\log|b_{i,1}\lambda_1+b_{i,2}\lambda_2|)
\end{align*}
So we can compute derivatives with respect to $\lambda_1$:
\begin{align}
\label{Eq::Diff of c1}
   \frac{\partial c_1}{\partial \lambda_1}  =\frac{(1-\sum\limits_{i=1}^{n}\beta_{i})^2}{(1-\sum\limits_{i=1}^{n}\beta_{i})\lambda_1+(1-\sum\limits_{i=1}^{n}\gamma_{i})\lambda_2} +\sum_{i=1}^{n}\frac{\beta_i^2}{\beta_i\lambda_1+\gamma_i\lambda_2}-\frac{1}{\lambda_1}
 \end{align}
 \begin{align}
\label{Eq::Diff of c2}
   \frac{\partial c_2}{\partial \lambda_1}= \frac{1-\sum\limits_{i=1}^{n}\gamma_{i}}{\lambda_1+\frac{1-\sum\limits_{i=1}^{n}\gamma_{i}}{1-\sum\limits_{i=1}^{n}\beta_{i}}\lambda_2} +\sum_{i=1}^{n}\frac{\gamma_i}{\lambda_1+\frac{\gamma_i}{\beta_i}\lambda_2}
 \end{align}

Now note that $\Gamma_\varepsilon(\mathcal{A},B)$ is a curve in $\mathbb{R}^2$. We say $\Gamma_\varepsilon(\mathcal{A},B)$ has a cusp at $\xi_{B,\varepsilon}(\lambda_1, \lambda_2)$ iff $\frac{\partial c_1}{\partial \lambda_1}(\lambda_1,\lambda_2)=\frac{\partial c_2}{\partial \lambda_1}(\lambda_1,\lambda_2)=0$.
However, in our setting, it the vanishing of just one partial derivative
suffices.
\begin{lem}
\label{Lem::Cusps occur}
$\Gamma_\varepsilon(\mathcal{A},B)$ has a cusp at $\xi_{B,\varepsilon}(\lambda_1, \lambda_2)$ iff $\frac{\partial c_1}{\partial \lambda_1}\frac{\partial c_2}{\partial \lambda_1}\!=\!0$. Moreover,
$\Gamma_\varepsilon(\mathcal{A},B)$ has at most $n$ cusps.
\end{lem}

\noindent
{\bf Proof:}
The first statement follows from Lemma \ref{Lem::Normal Vectors}: Since
$\lambda_1\frac{\partial c_1}{\partial \lambda_1}+\lambda_2\frac{\partial c_2}{\partial \lambda_1}=0$, one partial vanishing implies the other vanishes.
For the second statement, just observe that the numerator of $\displaystyle{\frac{\partial c_2}{\partial \lambda_1}}$ is a homogeneous polynomial with degree $n$ after clearing
denominators. \qed

\subsection{Reduced Chamber Structure}
\label{Chambers Structure}
Following \cite{RojasRusek_Adiscriminant}, we call the bounded (resp.\
unbounded) pieces of $\mathbb{R}^{k-1}\setminus
\Gamma_{\varepsilon}(\mathcal{A},B)$ \emph{signed reduced inner (resp.\ outer)
chambers}.
\begin{lem}
\cite[Proposition 4.6]{DengRojasTelek}
\label{Lem::One Cusp}
If $\cA$ is a near-circuit and
$\Gamma_{\varepsilon}(\mathcal{A},B)$ has at most one cusp
then the complement of the signed reduced $\mathcal{A}$-discriminant $\Gamma_\varepsilon(\mathcal{A},B)$ has at most two pieces, both of which are unbounded.
\end{lem}

We want to find the chamber structure when
$\Gamma_{\varepsilon}(\mathcal{A},B)$ has more cusps. Let us first recall
Cauchy's Mean Value Theorem.

\begin{lem}(Cauchy's Mean Value Theorem)
\label{lem:CMVT}
If the functions $f, g \colon [a,b] \to \mathbb{R}$ are both continuous and differentiable on the open interval $(a,b)$, then there exists some $c\in (a,b)$, such that
\[
(f(b)-f(a))g'(c)=(g(b)-g(a))f'(c).
\]
\end{lem}

The $\mathcal{A}$-discriminant curve $\Gamma_{\varepsilon}(\mathcal{A},B)$ is
piece-wise smooth and we call the each smooth part a \emph{curve segment}.
\begin{lem}
    \label{Lem::Intersections of Segments}
    Every pair of curve segments intersects at most once.
\end{lem}
\begin{proof}
    Suppose there are two curve segments intersect at two points $p_1$ and $p_2$. Then by Lemma \ref{lem:CMVT}, in each curve segment, there exists a point
such that the tangent line at this point is parallel to the line connecting $p_1$ and $p_2$. Therefore, we find two points in $\Gamma_{\varepsilon}(\mathcal{A},B)$ such that the tangent lines at these two points are parallel, which is contradictory to Lemma \ref{Lem::Normal Vectors}.
\end{proof}

For a point in a given chamber, consider the path from this point to a point in an outer chamber, and we suppose this path intersects with $\Gamma_{\varepsilon}(\mathcal{A},B)$ transversally. We define the \emph{depth} of this chamber as
the minimum number of intersections of the path and $\Gamma_{\varepsilon}(\mathcal{A},B)$. For example, the depth of an outer chamber is 0, and the depth of
an inner chamber adjacent to an outer chamber is 1.

\begin{lem}
\label{Lem::depth by cusps}
    Suppose $\Gamma_{\varepsilon}(\mathcal{A},B)$ has exactly $m$ cusps. Then
each chamber has depth no more than $\floor{m/2}$.
\end{lem}
\begin{proof}
We have already built $B$ so that it has bottom row $[0,-1]$.
Suppose $\Gamma_{\varepsilon}(\mathcal{A},B)$ is parametrized by $\bar{\xi}_{\varepsilon,B}(t)$ where
$t\in (t_0,t_\infty)$ and $\mathcal{C}$ is a chamber of $\mathbb{R}^2\setminus \Gamma_{\varepsilon}(\mathcal{A},B)$. Let $s$ be a curve segment of $\Gamma_{\varepsilon}(\mathcal{A},B)$ such that $s=\bar{\xi}_{\varepsilon,B}((t_1,t_2))$, where $(t_1,t_2)\subset (t_0,t_\infty)$. Given  $\delta\in \mathbb{R}$, we define the following path:
    \begin{align}\label{Eq::Segment extension}
   \gamma_{\delta,s}(t)= \left\{ \begin{array}{lc}
 \displaystyle{\frac{d \bar{\xi}_{\varepsilon,B}}{d t}(t_1)\cdot t + \bar{\xi}_{\varepsilon,B}(t_1)+(0,\delta)} & t_0 < t \leq t_1\\
  \displaystyle{ \bar{\xi}_{\varepsilon,B}(t)+(0,\delta)} & t_1 < t < t_2\\
  \displaystyle{\frac{d \bar{\xi}_{\varepsilon,B}}{d t}(t_2)\cdot t + \bar{\xi}_{\varepsilon,B}(t_2)+(0,\delta)} & t_2 \leq t < t_\infty
\end{array} \right.
    \end{align}
Then $(\gamma_{\delta,s}(t))_1$ and $(\gamma_{\delta,s}(t))_2$ are in $C^1(t_0,t_\infty)$ (i.e., their derivatives are continuous). We can choose a point $p\in \mathcal{C}$ and a real number $\delta$ such that $p\in\gamma_{\delta,s}(t)$ and $\gamma_{\delta,s}(t)$ intersects $\Gamma_{\varepsilon}(\mathcal{A},B)$ transversally. Let $p=\gamma_{\delta,s}(t_p)$. It suffices to show that the number of intersections between $\gamma_{\delta,s}(t)$ and $\Gamma_{\varepsilon}(\mathcal{A},B)$ is at most $m$. Consequently,  one of the path $\gamma_{\delta,s}((t_0,t_p))$ or $\gamma_{\delta,s}((t_p,t_\infty))$ must intersect $\Gamma_{\varepsilon}(\mathcal{A},B)$ at most $\displaystyle{\left\lfloor\frac{m}{2}\right\rfloor}$ times.

In fact, when $\delta\neq 0$, we have $\gamma_{\delta,s} \bigcap s=\emptyset$. Suppose, to the contrary, there exists $t_a\in (t_1,t_2)$ and $t_b\in (t_0,t_\infty)$ such that $\bar{\xi}_{\varepsilon,B}(t_a)=\gamma_{\delta,s}(t_b)$. Then we have
\[
\gamma_{\delta,s}(t_a)-\gamma_{\delta,s}(t_b)=(0,\delta)
\]
By Lemma \ref{lem:CMVT}, there exists $t_c\in (t_a,t_b)$ such that the tangent vector at $\gamma_{\delta,s}(t_c)$ is parallel to $(0,\delta)$. In other words, the normal vector at $\gamma_{\delta,s}(t_c)$ is $(1,0)$. However, by Lemma \ref{Lem::Normal Vectors}, the normal vector at $\gamma_{\delta,s}(t_c)$ is $(t_c,1)$, leading to a contradiction.

For the remaining $m$ curve segments of $\Gamma_{\varepsilon}(\mathcal{A},B)$, we show that $\gamma_{\delta,s}$ intersects each of them at most once. To see this, let us argue by contradiction: suppose $\gamma_{\delta,s}$ intersects a given segment at two distinct points $p_1$ and $p_2$. By Lemma \ref{lem:CMVT},  there exist a point on the curve segment and a point on $\gamma_{\delta,s}$ where the tangent lines are parallel to the line joining $p_1$ and $p_2$. However, there always exists a point in $s$ where the tangent vector is identical to that of $\gamma_{\delta,s}$. This leads to a situation where two points on $\Gamma_{\varepsilon}(\mathcal{A},B)$ have the same tangent vectors, which yields a contradiction.
\end{proof}

\begin{cor}
    \label{Lem::Structure of 2 or 3 cusps}
    If $\Gamma_{\varepsilon}(\mathcal{A},B)$ has exactly 2 or 3 cusps, then
all the inner chambers have depth 1.
\end{cor}

\begin{lem}
    \label{Lem::Structure of 4 cusps}
    If $\Gamma_{\varepsilon}(\mathcal{A},B)$ has exactly 4 cusps, then the depth of any chamber is at most $2$. Moreover, any chamber adjacent to a
depth-$2$ chamber must have depth 1.
\end{lem}
\begin{proof}
    The first statement follows directly from Lemma \ref{Lem::depth by cusps}. For the second statement, observe that any chamber adjacent to a depth-$2$ chamber must have depth at least $1$; otherwise, it would be an outer chamber, which cannot be adjacent to a depth-$2$ inner chamber.

    Now, suppose there exist two adjacent depth-$2$ inner chambers, denoted by $\mathcal{C}_1$ and $\mathcal{C}_2$, with a common boundary lying on a curve segment $s$. Let $p$ be a point on this common boundary, and choose a curve segment $s'$ distinct from $s$. Consider the curve $\gamma_{\delta, s'}(t)$ defined in \eqref{Eq::Segment extension}, chosen so that $p \in \gamma_{\delta, s'}(t)$. By the proof of Lemma \ref{Lem::depth by cusps}, $\gamma_{\delta, s'}(t)$ intersects $\Gamma_{\varepsilon}(\mathcal{A}, B)$ at most $4$ times, dividing $\gamma_{\delta, s'}(t)$ into $5$ segments. It follows that only one of these segments can have depth $2$. However, $\mathcal{C}_1 \cap \gamma_{\delta, s'}$ and $\mathcal{C}_2 \cap \gamma_{\delta, s'}$ belong to different segments and both have depth $2$, which leads to a contradiction.
\end{proof}

\section{Morse Theory and Hessians Along $\mathcal{A}$-discriminants}
\label{Morse}
In this section, we discuss the number of pieces of $Z_\R(f_c)$ for $c$ in
each discriminant chamber. Theorem \ref{Thm::RojasRusek_Adiscriminant} establishes that if two exponential sums belong to the same chamber, then their zero sets contain the same number of pieces. In Proposition \ref{lem::Cases when the number changes}, we describe how the number of pieces varies between adjacent chambers.

Morse Theory plays a crucial role in proving these results. We begin by recalling the statements of the Morse Lemma and Morse Theorem (see \cite{Book::MorseTheory}).

\begin{lem}(Morse Lemma)
\label{Lem::MorseLemma}
    Let $p$ be a non-degenerate critical point for a Morse function $f$. Then there is a local coordinate system $(y_1,\ldots,y_n)$ in a neighborhood $U$ of $p$ with $y_i(p)=0$ for all $i$ and such that the identity
    \[
    f(y)=f(p)-y_1^2+\cdots-y_s^2+y_{s+1}^2+\cdots+y_n^2
    \]
    holds throughout $U$, where $s$ is the number of negative eigenvalues of Hessian of $f$ at $p$. We call $s$ the index of $f$ at $p$
\end{lem}

By Lemma \ref{Thm::RojasRusek_Adiscriminant}, it is shown that if two exponential sums lie in the same chamber, then they share the same isotopy type, as established through Morse Theory. In particular, this implies that they have the same number of pieces. Thus, it is meaningful to refer to the number of pieces in each chamber.

Now let us focus on when the numbers change for two adjacent chambers. Let
$f_c$ be an exponential sum with $\alpha_1=0$ in its support. Then on
$Z_{\mathbb{R}}(f_c)$ we always have $-c_1=\sum_{i=2}^{n+k}
c_ie^{\alpha_i\cdot x}$.
By varying the value of $c_1$ while keeping $c_2,\ldots,c_{n+k}$ fixed, we trace out a line in the reduced ambient space $\mathbb{R}^{k-1}$. If this line intersects $\Gamma_{\varepsilon}(\mathcal{A},B)$ transversally, the corresponding exponential sums transition from one chamber to another. In such cases, we can examine the change in pieces by applying the Morse Lemma at the intersection.

\begin{prop}
\label{lem::Cases when the number changes}
Let $f_c$ be the exponential sums as defined before. Assume that the line obtained by varying $c_1$ intersects $\Gamma_{\varepsilon}(\mathcal{A},B)$ transversely at $f_{c^*}$. For the two adjacent chambers with $f_{c^*}$ on their shared boundary,
    the number of pieces changes only when the index (i.e., the number of
negative eigenvalues of the Hessian) of $f_{c^*}$ at the critical point
$x^*$ is one of the following: $0$, $1$, $n-1$, or $n$.
\end{prop}

\noindent
{\bf Proof:} First, note that $\frac{\partial f_c}{\partial x_i} =
\frac{\partial f_{c^*}}{\partial x_i} = -\frac{\partial c_1}{\partial x_i} $
for all \( i \). Consider the Morse function that maps the zero set
\[ \{(x_1, \ldots, x_n, c_1) \in \mathbb{R}^{n+1} \mid f_c = 0\} \]
to \( c_1 \). By the Morse Lemma (Lemma \ref{Lem::MorseLemma}), in a small
neighborhood of the critical point \( (x^*, c_1^*) \), we can find a new chart
\( (y_1, \ldots, y_n) \) such that
$c_1 - c_1^* = -y_1^2 - y_2^2 - \cdots - y_s^2+ y_{s+1}^2+\cdots + y_n^2$,
where  \( s \) is the index of $f_{c^*}$.

In this quadratic form, as \( c_1 \) increases across \( c_1^* \), the isotopy type of the set of \( (y_1, \ldots, y_n) \) changes as follows:
\begin{itemize}
    \item If \( s = 0 \) or \( s = n \), the set changes from an empty set to an \( (n-1) \)-sphere (or vice versa).
    \item If \( s = 1 \) or \( s = n-1 \), the set changes from a hyperboloid of two sheets to a hyperboloid of one sheet (or vice versa). Both of these cases change the number of pieces.
    \item However, if \( 2 \leq s \leq n-2 \), the isotopy type of the set of \( (y_1, \ldots, y_n) \) changes from \( \mathbb{R}^{n-s} \times \mathbb{S}^{k-1} \) to \( \mathbb{R}^s \times \mathbb{S}^{n-s-1} \), which does not change the number of pieces. \qed
\end{itemize}

Thanks to Proposition \ref{lem::Cases when the number changes}, we can
determine the number of pieces of $Z_\R(f_c)$ for $c$ in each chamber by
computing the index of the exponential sums that correspond to a point on the \( \mathcal{A} \)-discriminant. From this point onward, we focus on honest
$n$-variate exponential sums with $(n+3)$ terms and non-degenerate supports.
As noted earlier, we can reduce any
near-circuit exponential sum to a special form so that
the Gale dual matrix $B$ has the form shown in
\eqref{Eq::ReducedGualDualMatrix}. Let us now compute the Hessian of $f_c$ when
$f_c$ has a singular zero.
\begin{lem}
\label{Lem::Critical point and Hessian} Let $f$ and $B$ be defined as above.
If  $(c_1,c_2)\in\Gamma_{\varepsilon}(\mathcal{A},B)$ with
$ (c_1,c_2)=\Log|(\lambda_1,\lambda_2)\cdot B^{\top}|\cdot B$, then $f$ has a
critical point $x^*$ with $x_i^*=\log |\beta_i \lambda_1+\gamma_i\lambda_2|-\log|(1-\sum_{i=1}^n \beta_i)\lambda_1+(1-\sum_{i=1}^n \gamma_i)\lambda_2|$ and
$f(x^*)=0$. The Hessian of $f$ at $x^*$ is given by
    \[
    \frac{1}{(1-\sum_{i=1}^n \beta_i)\lambda_1+(1-\sum_{i=1}^n \gamma_i)\lambda_2}(M(\beta)\lambda_1+M(\gamma)\lambda_2),
    \]
where $M(\beta):=\begin{bmatrix} \beta_1(\beta_1-1) & \beta_1\beta_2 & \cdots & \beta_1\beta_{n}\\
\beta_1\beta_2 & \beta_2(\beta_2-1) & \cdots & \beta_2\beta_{n}\\
\vdots & \vdots & \ddots & \vdots\\
\beta_1\beta_{n} & \beta_2\beta_{n} & \cdots & \beta_{n}(\beta_{n}-1)
\end{bmatrix}$
\end{lem}
\begin{proof}
    It is not hard to check that
    \[
    f(x^*)=\frac{\partial f}{\partial x_1}(x^*)=\cdots=\frac{\partial f}{\partial x_n}(x^*)=0
    \]
    by direct computation. For the Hessian, note that for all $i$ we have
$\varepsilon_ie^{x_i^*}=-\beta_i\varepsilon_{n+1}e^{\beta\cdot x^*-c_1}-\gamma_i\varepsilon_{n+2}e^{\gamma\cdot x^*-c_2}$
since $\frac{\partial f}{\partial x_i}(x^*)=0$.
\end{proof}

As derived earlier, we can assume \( \lambda_1 = \mu \) and
\( \lambda_2 = 1 \) since $[0, -1]$ is a row vector of \( B \). Let us now find
the characteristic polynomial of the Hessian obtained in Lemma \ref{Lem::Critical point and Hessian}.

\begin{lem}
\label{lem::Char_Poly}
The characteristic polynomial of $(M(\beta)\mu+M(\gamma))$ with $\mu\neq 0$, is
$p_\mu(\zeta)=\mu\prod_{i=1}^{n}(\zeta+\beta_i \mu +\gamma_i)g(\zeta, \mu)$,
where $g(\zeta,\mu)$ is the determinant of \\
$\begin{bmatrix}
\sum\limits_{i=1}^n\frac{\beta_i^2 }{\zeta+\beta_i\mu+\gamma_i}-\frac{1}{\mu} & -(1-\sum\limits_{i=1}^n\beta_i)-\zeta\sum\limits_{i=1}^n \frac{\beta_i}{\beta_i\mu + \gamma_i +\zeta}\\   -(1-\sum\limits_{i=1}^n\beta_i)-\zeta\sum\limits_{i=1}^n \frac{\beta_i}{\beta_i\mu + \gamma_i +\zeta} & -(1-\sum\limits_{i=1}^n\beta_i)\mu-(1-\sum\limits_{i=1}^n\gamma_i)-\zeta\sum\limits_{i=1}^n \frac{\beta_i\mu+\gamma_i}{\beta_i\mu + \gamma_i +\zeta}
\end{bmatrix}$.
\end{lem}

\noindent
{\bf Proof:} First note that
$M(\beta)=\beta^{\top}\beta-\diag(\beta_1,\ldots,\beta_n)$. Then
    \begin{align*}
        p_\mu(\zeta)&=\det(\zeta I_n - M(\beta)\mu-M(\gamma))\\
        &=\det(\diag(\zeta+\beta_1 \mu +\gamma_1,\ldots,\zeta+\beta_n \mu +\gamma_n)-(\beta^{\top}\beta\mu+\gamma^{\top}\gamma))\\
        &=\prod_{i=1}^{n}(\zeta+\beta_i \mu +\gamma_i)\cdot\\
        &\hspace{-13mm}\det\left(I_n-\diag\left((\zeta+\beta_1 \mu +\gamma_1)^{-1}, \ldots, (\zeta+\beta_n \mu +\gamma_n)^{-1}\right)
        \begin{bmatrix}
            \beta^{\top} & \gamma^{\top}
        \end{bmatrix}
        \begin{bmatrix}
            \beta \mu\\
            \gamma
        \end{bmatrix}\right)\\
        &=\prod_{i=1}^{n}(\zeta+\beta_i \mu +\gamma_i)\cdot\\
        &\hspace{-13mm}\det\left(I_2-\begin{bmatrix}
            \beta \mu\\
            \gamma
        \end{bmatrix}
        \diag\left((\zeta+\beta_1 \mu +\gamma_1)^{-1}, \ldots, (\zeta+\beta_n \mu+\gamma_n)^{-1}\right)
        \begin{bmatrix}
            \beta^{\top} & \gamma^{\top}
        \end{bmatrix}\right)\\
        &=\mu\prod_{i=1}^{n}(\zeta+\beta_i \mu +\gamma_i)\det
        \begin{bmatrix}
         \sum\limits_{i=1}^n\frac{\beta_i^2 }{\zeta+\beta_i\mu+\gamma_i}-\frac{1}{\mu} &\sum\limits_{i=1}^n\frac{\beta_i\gamma_i}{\zeta+\beta_i\mu+\gamma_i}\\   \sum\limits_{i=1}^n\frac{\beta_i\gamma_i }{\zeta+\beta_i\mu+\gamma_i} &  \sum\limits_{i=1}^n\frac{\gamma_i^2}{\zeta+\beta_i\mu+\gamma_i} -1
        \end{bmatrix}
    \end{align*}
    Here we use Sylvester's determinant identity, which states that $\det(I_m-AB)=\det(I_n-BA)$  if $A$ and $B$
 are matrices of sizes $m\times n$
 and $n\times m$, respectively. The first expression of $g(\zeta,\mu)$ has already been derived.

Notice that we have the identities\\
$(1-\sum\limits_{i=1}^n\beta_i)+\mu\left(\sum\limits_{i=1}^n \frac{\beta_i^2}{\beta_i\mu + \gamma_i +\zeta} -\frac{1}{\mu}\right)+ \sum\limits_{i=1}^n \frac{\beta_i\gamma_i}{\beta_i\mu + \gamma_i +\zeta}+\zeta\sum\limits_{i=1}^n \frac{\beta_i}{\beta_i\mu + \gamma_i +\zeta}=0$
and\\
$(1-\sum\limits_{i=1}^n\gamma_i)+\mu\left(\sum\limits_{i=1}^n \frac{\beta_i\gamma}{\beta_i\mu + \gamma_i +\zeta} \right)+ \sum\limits_{i=1}^n \frac{\gamma_i^2}{\beta_i\mu + \gamma_i +\zeta}-1+\zeta\sum\limits_{i=1}^n \frac{\gamma_i}{\beta_i\mu + \gamma_i +\zeta}=0$.

Therefore,
 \begin{align*}
    & g(\zeta,\mu)\\
        &=\det
        \begin{bmatrix}
         \sum\limits_{i=1}^n\frac{\beta_i^2 }{\zeta+\beta_i\mu+\gamma_i}-\frac{1}{\mu} &\mu\left(\sum\limits_{i=1}^n\frac{\beta_i^2 }{\zeta+\beta_i\mu+\gamma_i}-\frac{1}{\mu}\right)+\sum\limits_{i=1}^n\frac{\beta_i\gamma_i}{\zeta+\beta_i\mu+\gamma_i}\\   \sum\limits_{i=1}^n\frac{\beta_i\gamma_i }{\zeta+\beta_i\mu+\gamma_i} &  \mu\left(\sum\limits_{i=1}^n\frac{\beta_i\gamma_i }{\zeta+\beta_i\mu+\gamma_i}\right)+\sum\limits_{i=1}^n\frac{\gamma_i^2}{\zeta+\beta_i\mu+\gamma_i} -1
        \end{bmatrix}\\
         &=\det
        \begin{bmatrix}
         \sum\limits_{i=1}^n\frac{\beta_i^2 }{\zeta+\beta_i\mu+\gamma_i}-\frac{1}{\mu} & -(1-\sum\limits_{i=1}^n\beta_i)-\zeta\sum\limits_{i=1}^n \frac{\beta_i}{\beta_i\mu + \gamma_i +\zeta}\\   \sum\limits_{i=1}^n\frac{\beta_i\gamma_i }{\zeta+\beta_i\mu+\gamma_i} &  -(1-\sum\limits_{i=1}^n\gamma_i)-\zeta\sum\limits_{i=1}^n \frac{\gamma_i}{\beta_i\mu + \gamma_i +\zeta}
        \end{bmatrix}\\
        &=\det
\text{$\begin{bmatrix}
         \sum\limits_{i=1}^n\frac{\beta_i^2 }{\zeta+\beta_i\mu+\gamma_i}-\frac{1}{\mu} & -(1-\sum\limits_{i=1}^n\beta_i)-\zeta\sum\limits_{i=1}^n \frac{\beta_i}{\beta_i\mu + \gamma_i +\zeta}\\   -(1-\sum\limits_{i=1}^n\beta_i)-\zeta\sum\limits_{i=1}^n \frac{\beta_i}{\beta_i\mu + \gamma_i +\zeta} & -(1-\sum\limits_{i=1}^n\beta_i)\mu-(1-\sum\limits_{i=1}^n\gamma_i)-\zeta\sum\limits_{i=1}^n \frac{\beta_i\mu+\gamma_i}{\beta_i\mu + \gamma_i +\zeta}
        \end{bmatrix}$}  \text{ \qed }
 \end{align*}

For each $\mu\in\mathbb{R}$, the parametric map $\bar{\xi}_{B,\varepsilon}$ provides a point on the $\mathcal{A}$-discriminant curve, which corresponds to an exponential sum, allowing us to compute its Hessian. Although we have derived the characteristic polynomial of the Hessian, determining the index directly remain a challenge. However, since $\mu$ varies along the $\mathcal{A}$-discriminant curve $\Gamma(\mathcal{A},B)$, the index changes only at the cusps.

\begin{lem}
\label{lem::Cusps and Hessian}
    For $\mu\neq 0$, the index of the Hessian $(M(\beta)\mu+M(\gamma))$ changes only at the cusps of $\mathcal{A}$-discriminant curve, where $\mu$ satisfies the condition $\displaystyle{\frac{\partial c_1}{\partial \lambda_1}(\mu,1)}=0$
\end{lem}

\noindent
{\bf Proof:} The sign of the eigenvalues changes only when the matrix becomes
singular, i.e., when $0$ is an eigenvalue. From the characteristic polynomial
(Lemma \ref{lem::Char_Poly}), we deduce that $p_\mu(0)$ is
$\mu\prod_{i=1}^{n}(\beta_i \mu +\gamma_i)\det
        \begin{bmatrix}
         \sum\limits_{i=1}^n\frac{\beta_i^2 }{\zeta+\beta_i\mu+\gamma_i}-\frac{1}{\mu} & -(1-\sum\limits_{i=1}^n\beta_i)\\   -(1-\sum\limits_{i=1}^n\beta_i) & -(1-\sum\limits_{i=1}^n\beta_i)\mu-(1-\sum\limits_{i=1}^n\gamma_i)
         \end{bmatrix}$.
Since $\mu\neq 0$, and $\beta_i\mu+\gamma_i\neq 0$ for all $i$, we have
$p_\mu(0)=0$ iff\\  $
\left(\sum\limits_{i=1}^n\frac{\beta_i^2 }{\zeta+\beta_i\mu+\gamma_i}-\frac{1}{\mu}\right)\left( -(1-\sum\limits_{i=1}^n\beta_i)\mu-(1-\sum\limits_{i=1}^n\gamma_i)\right)-(1-\sum\limits_{i=1}^n\beta_i)^2=0$. \\
That is,
$\sum\limits_{i=1}^n\frac{\beta_i^2 }{\zeta+\beta_i\mu+\gamma_i}-\frac{1}{\mu}+\frac{(1-\sum\limits_{i=1}^n\beta_i)^2}{(1-\sum\limits_{i=1}^n\beta_i)\mu+(1-\sum\limits_{i=1}^n\gamma_i)}=0$,
which corresponds exactly to $\displaystyle{\frac{\partial c_1}{\partial \lambda_1}(\mu,1)=0}$. Therefore, this condition defines the cusps of the $\mathcal{A}$-discriminant curve. \qed

\smallskip
We also show that $0$ is a simple root of the characteristic polynomial $p_\mu(\zeta)$ when $\mu$ satisfies $\displaystyle{\frac{\partial c_1}{\partial \lambda_1}(\mu,1)}=0$. This implies that at the cusps, at most one eigenvalue changes its sign.

\begin{lem}
\label{lem::zero single eigen}
    Let $\mu_0\neq 0$ be a root of $\displaystyle{\frac{\partial c_1}{\partial \lambda_1}(\mu,1)}$ with $\beta_i\mu_0+\gamma_i\neq 0$ for all $i$ and $\displaystyle{(1-\sum\limits_{i=1}^n\beta_i)\mu_0+(1-\sum\limits_{i=1}^n\gamma_i)\neq 0}$. Then $\displaystyle{\frac{\partial p_{\mu_0}}{\partial \zeta}(0)\neq 0}$.
\end{lem}
\begin{proof}
    $\zeta=0$ is a root of $p_{\mu_0}(\zeta)$ if and only if $(0,\mu_0)$ is a root of $g(\zeta,\mu)$ by our assumption. Also, we have
    \[
    \frac{\partial p_{\mu_0}}{\partial \zeta}(0)=\mu_0\prod_{i=1}^{n}(\zeta+\beta_i \mu_0 +\gamma_i)\frac{\partial g}{\partial \zeta}(0,\mu_0)
    \]
    By direct computation, one can find the derivative of $g$ with respect to $\zeta$ at the point $(0,\mu_0)$:
    \begin{align*}
    \frac{\partial g}{\partial \zeta}(0,\mu_0)&=\left((1-\sum\limits_{i=1}^n\beta_i)\mu_0+(1-\sum\limits_{i=1}^n\gamma_i)\right) \cdot
    \sum\limits_{i=1}^n\left(\frac{\beta_i}{\beta_i\mu_0 + \gamma_i }-\frac{1-\sum\limits_{j=1}^n\beta_j}{(1-\sum\limits_{j=1}^n\beta_j)\mu_0+(1-\sum\limits_{j=1}^n\gamma_j)}\right)^2
\end{align*}
By our assumption,
$\displaystyle{\frac{\partial g}{\partial \zeta}(0,\mu_0)}=0$ only when
\[\frac{\beta_i}{\beta_i\mu_0 + \gamma_i }-\frac{1-\sum\limits_{j=1}^n\beta_j}{(1-\sum\limits_{j=1}^n\beta_j)\mu_0+(1-\sum\limits_{j=1}^n\gamma_j)}=0\]
for all $i$, but this can't happen: If so, we would obtain
$\frac{\beta_1}{\gamma_1}=\frac{\beta_2}{\gamma_2}=\cdots
=\frac{\beta_n}{\gamma_n}$, which contradicts our setting of $\beta$ and
$\gamma$.
\end{proof}

In Lemma \ref{lem::zero single eigen}, we proved that only one eigenvalue may
change its sign at a cusp. Setting
$p_\mu(\zeta)\!=\!0$ then defines an implicit function $\zeta(\mu)$, with
$\zeta(\mu_0)=0$, that is well-defined in a neighborhood of $(0,\mu_0)$ by
Lemma \ref{lem::zero single eigen}.
Moreover, we can detect the sign change of this eigenvalue by analyzing the derivatives of $c_1$ with respect to $\lambda_1$. The following lemma provides further details:

\begin{lem}
\label{lem::Sign Changes}
 Let $\mu_0\neq 0$ be a root of $\displaystyle{\frac{\partial c_1}{\partial \lambda_1}(\mu,1)}$ as stated in Lemma \ref{lem::zero single eigen}, and let $\zeta(\mu)$ be the implicit function defined above. Then, for a sufficiently small neighborhood $\mu\in (\mu_0-\delta,\mu_0+\delta)$ of $\mu_0$, we have
 \[
    \sign(\zeta(\mu))=\sign\left(\frac{\partial c_1}{\partial \lambda_1}(\mu,1)\right).
    \]

\end{lem}

\noindent
{\bf Proof:} Suppose $\mu_0$ is a root of $\displaystyle{\frac{\partial c_1}
{\partial \lambda_1}(\mu,1)}$ with multiplicities $l$, which yields
$\frac{\partial^{l+1} c_1}{\partial \lambda_1^{l+1}}(\mu_0,1)\neq 0$
and $\frac{\partial^j c_1}{\partial \lambda_1^j}(\mu_0,1)=0$
for all $1\leq j\leq l$.

    Let $p(\zeta,\lambda)=\lambda\prod_{i=1}^{n}(\zeta+\beta_i \lambda +\gamma_i)g(\zeta,\lambda)$, where $g$ is the same as in Lemma
\ref{lem::Cusps and Hessian}.
   Also, by the proof of Lemma \ref{lem::Cusps and Hessian}, one can show that
   \[
   g(0,\mu)=-\left((1-\sum\limits_{i=1}^n\beta_i)\mu+(1-\sum\limits_{i=1}^n\gamma_i)\right)\frac{\partial c_1}{\partial \lambda_1}(\mu,1).
   \]
   Then we have
    \begin{align*}
        \frac{\partial^l g}{\partial \mu^l}(0,\mu_0)=
        -\left((1-\sum\limits_{i=1}^n\beta_i)\mu_0+(1-\sum\limits_{i=1}^n\gamma_i)\right)\frac{\partial^{l+1} c_1}{\partial \lambda_1^{l+1}}(\mu_0,1)\quad
        \end{align*}
        and
$\frac{\partial^j g}{\partial \mu^j}(\mu_0,1)=0$
for all $0\leq j\leq l-1$.

Therefore, on a small neighborhood of $\mu_0$, we have the following
(by induction on the order of the derivatives):
\[
\frac{\partial^j \zeta}{\partial \mu^j}(\mu_0)=-\frac{\frac{\partial^j p}{\partial \mu^j}(0,\mu_0)}{\frac{\partial p}{\partial \zeta}(0,\mu_0)}=-\frac{\frac{\partial^j g}{\partial \mu^j}(0,\mu_0)}{\frac{\partial g}{\partial \zeta}(0,\mu_0)}=0
\]
for all $0\leq j\leq l-1$. Also,
\[
\frac{\partial^l \zeta}{\partial \mu^l}(\mu_0)=-\frac{\frac{\partial^l g}{\partial \mu^l}(0,\mu_0)}{\frac{\partial g}{\partial \zeta}(0,\mu_0)}=\frac{\frac{\partial^{l+1} c_1}{\partial \lambda_1^{l+1}}(\mu_0,1)}{\sum\limits_{i=1}^n\left(\frac{\beta_i}{\beta_i\mu_0 + \gamma_i }-\frac{1-\sum\limits_{j=1}^n\beta_j}{(1-\sum\limits_{j=1}^n\beta_j)\mu_0+(1-\sum\limits_{j=1}^n\gamma_j)}\right)^2}
\]

Now let us consider the Taylor expansions of $\zeta(\mu)$ and
$\frac{\partial c_1}{\partial \lambda_1}(\mu,1)$ centered at $\mu_0$ within a small neighborhood $(\mu_0-\delta,\mu_0+\delta)$. Based on the calculations above, the leading terms of these expansions are $\frac{\partial^l \zeta}{\partial \mu^l}(\mu_0)\frac{(\mu-\mu_0)^l}{l!}$ and  $\frac{\partial^{l+1} c_1}{\partial \lambda_1^{l+1}}(\mu_0,1)\frac{(\mu-\mu_0)^l}{l!}$ , respectively. These terms dominate the signs of $\zeta(\mu)$ and $\displaystyle{\frac{\partial c_1}{\partial \lambda_1}(\mu,1)}$ in the neighborhood $(\mu_0-\delta,\mu_0+\delta)$. Since $\displaystyle{\frac{\partial^l \zeta}{\partial \mu^l}(\mu_0)}$ and $\displaystyle{\frac{\partial^{l+1} c_1}{\partial \lambda_1^{l+1}}(\mu_0,1)}$ share the same sign within this neighborhood, the proof follows. \qed

\section{Reduced $\mathcal{A}$-discriminants with multiple cusps}\label{Section::multiple cusps}
We now study near-circuit exponential sums under the additional assumption that their signed reduced $\mathcal{A}$-discriminant curves has at least two cusps.
To count the pieces of their real zero sets, we first consider the case where the $\mathcal{A}$-discriminant curve contains at least five cusps. We show that
the signature of the Hessian is sufficiently refined, ensuring that all chambers have the same number of pieces. Subsequently, we examine the scenario where the $\mathcal{A}$-discriminant curve has at most four cusps, analyzing the pieces through the structure of the chambers.

We begin by recalling a useful result of Bihan, Humbert, and Tavenas.
\begin{lem}
\cite[Prop.\ 6.2]{BihanBound}
\label{Lem::Number in Outer Chambers}
Suppose $f_c$ is a near-circuit exponential sum, corresponding to a point in
an outer chamber. Then $Z_\R(f_c)$ has at most 2 pieces.
\end{lem}

By combining Lemma \ref{Lem::Number in Outer Chambers} and \ref{Lem::diff of adjacent chambers}, we obtain a rough bound on the number of pieces.

\begin{cor}
    \label{Cor::Rough bound}
    If $f_c$ lies in a reduced signed chamber of depth $d$ then $Z_\R(f_c)$
has at most $2+d$ pieces.
\end{cor}

Using Proposition \ref{lem::Cases when the number changes}, we can more precisely count the pieces in the inner chambers by determining the index of $f_c$ when $f_c$ lies on $\mathcal{A}$-discriminant curve. While it may be challenging to directly compute the index of $f_c$ at a general point on the $\mathcal{A}$-discriminant curve, Lemma \ref{lem::Cusps and Hessian} and Lemma \ref{lem::Sign Changes} allow us to determine the index of $f_c$ as $\mu\to 0$, and track the sign changes of eigenvalues as $\mu$ changes.

Before presenting the main result, we first recall the setting of the exponential sums. We use the canonical form of $f_c$ as given in \eqref{Eq::ReducedExpSum} and the Gale dual matrix in \eqref{Eq::ReducedGualDualMatrix}. If $f_c$ contains a singular zero, the Hessian of $f_c$ is $M(\beta)\mu+M(\gamma)$, up to a scalar, as shown in Lemma \ref{Lem::Critical point and Hessian}. By taking the limit $\mu\to 0$,  we obtain $M(\gamma)$, whose index can be determined.

\begin{lem}
\label{lem::Signature of Hessian}
    The signs of the eigenvalues of $M(\gamma)$ (counted with
multiplicities) are
$\sign(\sum\limits_{i=1}^n \gamma_i-1, -\gamma_1, \ldots, -\gamma_n)
\setminus\{-\}$.
\end{lem}

\noindent
{\bf Proof:}
    By the proof of Lemma \ref{lem::Char_Poly}, as $\mu \to 0$, the characteristic polynomial is given by
    \[
    p_0(\zeta)=\prod_{i=1}^{n}(\zeta+\gamma_i)\left(1-\sum\limits_{i=1}^n\frac{\gamma_i^2}{\zeta+\gamma_i}\right)
    \]

We first assume that the $\gamma_i$'s are distinct. In this case, if $\zeta=-\gamma_j$ for some $j$, then
\[
p_0(-\gamma_j)=\gamma_j^2\prod_{i\neq j}(\gamma_i-\gamma_j)
\]
Since $\gamma_j\neq 0$  (by the assumption of the Gale dual matrix $B$) and $\gamma_i\neq \gamma_j$, it follows that $-\gamma_j$ is not a root of $p_\mu(\zeta)$.

As $p_\mu(\zeta)$ is a polynomial of degree $n$, its $n$ roots must all come from solving $\bar{p}_0(\zeta):=1-\sum\limits_{i=1}^n\frac{\gamma_i^2}{\zeta+\gamma_i}=0$. Note that \(\bar{p}\) has \(n\) poles at
\(-\gamma_1, \ldots, -\gamma_n\). Without loss of generality, assume that
\[
-\gamma_1 < -\gamma_2 < \cdots < -\gamma_k < 0 < -\gamma_{k+1} < \cdots < -\gamma_n
\]
for some \(1 \leq k \leq n\). Moreover, we have
$\bar{p}'_0(\zeta) = \sum_{i=1}^n \frac{\gamma_i^2}{(\zeta + \gamma_i)^2} > 0$,
which shows that \(\bar{p}\) is strictly increasing in each interval between its poles.

At each pole \(-\gamma_j\), the behavior of \(\bar{p}\) is characterized as follows:
\[
\lim_{\zeta \to -\gamma_j^-} \bar{p}_0(\zeta) = +\infty, \quad \lim_{\zeta \to -\gamma_j^+} \bar{p}_0(\zeta) = -\infty.
\]
Thus, \(\bar{p}\) has the following roots:
\begin{itemize}
    \item \(k-1\) \textbf{negative roots}, respectively located in the intervals \[(- \gamma_1, - \gamma_2), \ldots, (- \gamma_{k-1}, - \gamma_k);\]
    \item \(n-k-2\) \textbf{positive roots}, respectively located in the intervals \[(- \gamma_{k+1}, - \gamma_{k+2}), \ldots, (- \gamma_{n-1}, - \gamma_n);\]
    \item One \textbf{positive root} in \((- \gamma_n, +\infty)\), since \(\lim_{\zeta \to \infty} \bar{p}_0(\zeta) = 1\).
\end{itemize}

The sign of the root in the interval \((- \gamma_k, - \gamma_{k+1})\) depends
on the value of \(\bar{p}_0(0)\): $\bar{p}_0(0) = 1 - \sum_{i=1}^n \gamma_i$.
\begin{itemize}
    \item If \(\bar{p}_0(0) > 0\), the root in this interval is negative.
    \item Otherwise, the root is positive.
\end{itemize}

Since \(k\) is the number of positive \(\gamma_i\)'s, we can summarize the signs of the roots as follows:  The signs of the roots are exactly
\[
\sign\left(\sum_{i=1}^n \gamma_i - 1, -\gamma_1, \ldots, -\gamma_n\right),
\]
after removing one negative sign.

Now we consider the case when the \(\gamma_i\)'s are not distinct. Suppose that \(\gamma_1, \ldots, \gamma_l\) are distinct with \(1 \leq l \leq n\), and
\[
p_0(\zeta) = \prod_{i=1}^l (\zeta + \gamma_i)^{k_i} \left( 1 - \sum_{i=1}^l \frac{k_i \gamma_i^2}{\zeta + \gamma_i} \right),
\]
where \(k_i \geq 1\) is the multiplicity of \(\gamma_i\).

Now, consider the polynomial
\[
\hat{p}_0(\zeta) := \frac{p_0(\zeta)}{\prod\limits_{i=1}^l (\zeta + \gamma_i)^{k_i-1}} = \prod_{i=1}^l (\zeta + \gamma_i) \left( 1 - \sum_{i=1}^l \frac{k_i \gamma_i^2}{\zeta + \gamma_i} \right).
\]
The signs of the roots of \(\hat{p}_0(\zeta)\) correspond to the case we discussed previously, where all \(\gamma_i\)'s ($1\leq i\leq l$) are distinct.

Additionally, \(p_\mu(\zeta)\) has roots \(-\gamma_i\) with multiplicities \(k_i - 1\). Therefore, the signs of the roots are still given by
\[
\sign\left(\sum_{i=1}^n \gamma_i - 1, -\gamma_1, \ldots, -\gamma_n \right),
\]
after removing one negative sign. \qed

\smallskip
Consider now the case where the signed reduced
\(\mathcal{A}\)-discriminant curve \(\Gamma_{\varepsilon}(\mathcal{A},B)\) of \(f_c\) is parameterized by \(\mu \in (0, \infty)\), i.e.,
\[
\Gamma_{\varepsilon}(\mathcal{A},B) = \{\xi_{\varepsilon,B}(\mu,1) \mid \mu \in (0, \infty)\},
\]
(the case when \(\mu \in (-\infty, 0)\) is similar).

\begin{lem}
\label{lem::Signs by Cusps}
Suppose the reduced \(\mathcal{A}\)-discriminant curve \(\Gamma(\mathcal{A},B)\) has no poles when \(\mu > 0\), and \(\Gamma_{\varepsilon}(\mathcal{A},B)\) is the signed reduced \(\mathcal{A}\)-discriminant curve for \(\mu \in (0, \infty)\). If \(\Gamma_{\varepsilon}(\mathcal{A},B)\) has \(m\) cusps (counted with
multiplicities, and \(2 \leq m \leq n\)), then
$\sign\left(1 - \sum_{i=1}^{n}\gamma_{i}, \gamma_1, \ldots, \gamma_n \right)$
contains \(\displaystyle{\left\lfloor\frac{m}{2}\right\rfloor + 1}\) positive signs and \(\displaystyle{\left\lfloor\frac{m+1}{2}\right\rfloor}\) negative signs.
\end{lem}

\begin{proof}
By Lemma \ref{Lem::Cusps occur}, the reduced $\mathcal{A}$-discriminant curve has a cusp if and only if $\frac{\partial c_2}{\partial \lambda_1} = 0$.
Let $\mu_i$ $(1 \leq i \leq m)$ be the positive roots of $\frac{\partial c_2}{\partial \lambda_1}$. Also, $\frac{\partial c_2}{\partial \lambda_1}$ is a rational function whose numerator has degree $n$, which yields that $\displaystyle{\frac{\partial c_2}{\partial \lambda_1}}$ has at most $n-m$ negative roots.

We call $\mu$ a pole of $\displaystyle{\frac{\partial c_2}{\partial \lambda_1}}$ if $b_{i,1}\mu+b_{i,2}$ for some $i$. Let $b_{ij}$ $(1 \leq i \leq n+3, \, j = 1, 2)$ be the entries of $B$. Then by our assumption, $\displaystyle{\frac{\partial c_2}{\partial \lambda_1}}$ has $n+1$ negative poles:
$-\frac{b_{i,2}}{b_{i,1}}$ for $i\!\in\!\{1,\ldots,n+1\}$.
These $n+1$ poles cut the negative axis $(-\infty, 0)$ into $(n+2)$
sub-intervals, or $n$ sub-intervals excluding those with endpoints $-\infty$ or
$0$.

\smallskip
\textbf{Case 1: $m$ is odd.} Since $\frac{\partial c_2}{\partial \lambda_1}$
has at most $n-m$ negative roots, there are at least $m$ sub-intervals without any real roots of $\frac{\partial c_2}{\partial \lambda_1}$. Among these $m$
sub-intervals, there exist $(m+1)/2$ of them that are not
adjacent when $m$ is odd.

Suppose the endpoints of these $(m+1)/2$ intervals are
\[
\left(-\frac{b_{i_1,2}}{b_{i_1,1}}, -\frac{b_{i_2,2}}{b_{i_2,1}}\right), \left(-\frac{b_{i_3,2}}{b_{i_3,1}}, -\frac{b_{i_4,2}}{b_{i_4,1}}\right), \ldots, \left(-\frac{b_{i_m,2}}{b_{i_m,1}}, -\frac{b_{i_{m+1},2}}{b_{i_{m+1},1}}\right).
\]
We claim that $\mathrm{sign}(b_{i_1,2}, b_{i_2,2}, \dots, b_{i_{m+1},2})$ contains $\displaystyle{\frac{m+1}{2}}$ positive signs and $\displaystyle{\frac{m+1}{2}}$ negative signs. In fact, on each interval
$\left(-\frac{b_{i_{2t-1},2}}{b_{i_{2t-1},1}},
-\frac{b_{i_{2t},2}}{b_{i_{2t},1}}\right)$ (for $t\!\in\!\{1,\ldots,(m+1)/2\}$),
we have that $\frac{\partial c_2}{\partial \lambda_1}$ has no root, and hence
the sign does not change. We now show that $b_{i_{2t-1},2}$ and $b_{i_{2t},2}$
have different signs: Assume without loss of generality that
$\frac{\partial c_2}{\partial \lambda_1} > 0$ on $\left(-\frac{b_{i_{2t-1},2}}{b_{i_{2t-1},1}}, -\frac{b_{i_{2t},2}}{b_{i_{2t},1}}\right)$.
Therefore,
\[
\lim_{\mu \to -\frac{b_{i_{2t-1},2}}{b_{i_{2t-1},1}}^+} \frac{\partial c_2}{\partial \lambda_1}(\mu, 1) =
\lim_{\mu \to -\frac{b_{i_{2t},2}}{b_{i_{2t},1}}^-} \frac{\partial c_2}{\partial \lambda_1}(\mu, 1) = +\infty.
\]

By \eqref{Eq::Diff of c2}, the expression of
$\frac{\partial c_2}{\partial \lambda_1}(\mu, 1)$ is a sum of fractions. So
the signs as $\mu \to -\frac{b_{i_{2t-1},2}}{b_{i_{2t-1},1}}^+$ (resp. $\mu \to -\frac{b_{i_{2t},2}}{b_{i_{2t},1}}^-$) only depend on the term
$\frac{b_{i_{2t-1},2}}{\mu + \frac{b_{i_{2t-1},2}}{b_{i_{2t-1},1}}}$
(resp.\ $\frac{b_{i_{2t},2}}{\mu + \frac{b_{i_{2t},2}}{b_{i_{2t},1}}}$).
Hence, $b_{i_{2t-1},2} > 0$ and $b_{i_{2t},2} < 0$. This concludes the proof.

\textbf{Case 2: $m$ is even.}
Let $(-b_{i_l,2}/b_{i_l,1},0)$ be the
rightmost sub-interval.  If there is a root of
$\frac{\partial c_2}{\partial \lambda_1}$ on
$(-b_{1,2}/b_{1,1},0)$,
then by the proof above, there are $m+1$ sub-intervals not containing any real
roots of $\frac{\partial c_2}{\partial \lambda_1}$, so there are
$(m+2)/2$ of them are not adjacent, which give
$(m+2)/2$ positive signs and $(m+2)/2$ negative signs as the proof above.

Now we just need to consider the case when $\frac{\partial c_2}{\partial
\lambda_1}$ has no roots on $(-b_{i_l,2}/b_{i_l,1},0)$. Note that by
\eqref{Eq::Diff of c2}, we have
$\frac{\partial c_2}{\partial \lambda_1}(0,1)=1>0$.
Hence $\frac{\partial c_2}{\partial \lambda_1}>0$ on
$(-b_{i_l,2}/b_{i_l,1},0)$. Thus,
\[
\lim_{\mu\to -\frac{b_{i_l,2}}{b_{i_l,1}}^+}\frac{\partial c_2}{\partial \lambda_1}(\mu,1)=+\infty
\]
and that means $b_{i_l,2}>0$. Similar to the prove above, precluding the
intervals with endpoints $-\infty$ or $0$, there are $m$ sub-intervals without
any real roots of $\frac{\partial c_2}{\partial \lambda_1}$. Among these $m$
sub-intervals, there are $m/2$ sub-intervals are not adjacent to each other and
each of them are not adjacent to $(-b_{i_l,2}/b_{i_l,1},0)$. Therefore, we have
$m/2$ negative $b_{i,2}$ and $m/2$ positive $b_{i,2}$. In addition to
$b_{i_l,2}>0$, there are $\frac{m}{2}+1$ positive entries and
$m/2$ negative entries in
$(b_{12},b_{22},\ldots,b_{n+1,2})=(1-\sum\limits_{i=1}^{n}\gamma_{i},\gamma_1,\ldots,\gamma_n)$.
\end{proof}

\begin{cor}
\label{cor::Sign of Eigenvalues}
    Under the same assumption as Lemma \ref{lem::Signs by Cusps}, if
$\Gamma_{\varepsilon}(\mathcal{A},B)$ has $m$ cusps (counted with multiplicities), then $M(\gamma)$ has at least $\floor{\frac{m+1}{2}}$ positive eigenvalues
and $\floor{\frac{m}{2}}$ negative eigenvalues.
\end{cor}
\begin{proof}
    This follows directly from Lemma \ref{lem::Signature of Hessian} and Lemma \ref{lem::Signs by Cusps}
\end{proof}

\noindent
{\bf Proof of Theorem \ref{Thm::Main}:}
By Lemma \ref{Lemma::Transform}, we can suppose $f$ has the form in
\eqref{Eq::ReducedExpSum}
\[
f=\varepsilon_0+\varepsilon_{1}e^{x_1}+\cdots+\varepsilon_{n}e^{x_n}+\varepsilon_{n+1}e^{\beta\cdot x-c_1}+\varepsilon_{n+2}e^{\gamma\cdot x-c_2}
\]
and with the Gale dual matrix $B$ in \eqref{Eq::ReducedGualDualMatrix}.
Additionally, the signed reduced $\mathcal{A}$-discriminant
$\Gamma_{\varepsilon}(\mathcal{A},B)$ of $f$ is given by the
parametric curve
$(c_1,c_2)=\Log|(\mu,1)\cdot B^{\top}|\cdot B$,
where $\mu$ is the parameter and $\mu\in (-\infty,0)$ or $\mu\in (0,+\infty)$.
Let us discuss the cases of the number of the cusps (counting multiplicities) in $\Gamma_{\varepsilon}(\mathcal{A},B)$.
\begin{itemize}[leftmargin=*]
    \item \textbf{$\Gamma_{\varepsilon}(\mathcal{A},B)$ has at least 6 cusps.} By Corollary \ref{cor::Sign of Eigenvalues}, $M(\gamma)$ has at least 3 positive eigenvalues and 3 negative eigenvalues. Since $\displaystyle{\frac{\partial c_1}{\partial \lambda_1}(\mu, 1)}$ is a continuous function on $(-\infty,0)$ or $(0,+\infty)$, by Lemma \ref{lem::Sign Changes}, there are at least 2 positive eigenvalues and 2 negative eigenvalues in the Hessian $M(\beta)\mu+M(\gamma)$ for any $\mu\in (-\infty,0)$ or $\mu\in (0,+\infty)$. Therefore, by Lemma \ref{lem::Cases when the number changes} and Lemma \ref{Lem::Number in Outer Chambers}, the real zero set of $f$ has at most 2 pieces.

    \item \textbf{$\Gamma_{\varepsilon}(\mathcal{A},B)$ has exactly 5 cusps.} We first show that this case only occurs when $\mu\in (-\infty,0)$. Indeed, if $\mu\in (0,+\infty)$, then $\displaystyle{\frac{\partial c_2}{\partial \lambda_1}(\mu, 1)}$ has exactly 5 positive roots. Since $\displaystyle{\frac{\partial c_2}{\partial \lambda_1}(\mu, 1)}$ is a rational function and all the poles are negative, then it has the form
    \[
    \frac{\partial c_2}{\partial \lambda_1}(\mu, 1)=\frac{(\mu-\mu_1)\cdots (\mu-\mu_5)(\mu+\mu_6)\cdots (\mu+\mu_\ell)}{(\mu+p_1)\cdots (\mu+p_{n+1})}h(\mu),
    \]
    where $\mu_i>0$ and $p_i>0$ for all $i$, and $h(\mu)$ is a monic polynomial without any real roots. Since $\displaystyle{\lim_{\mu \to +\infty}h(\mu)=+\infty}$, then $h(0)>0$. Thus we have $\frac{\partial c_2}{\partial \lambda_1}(0, 1)=-\frac{\mu_1\cdots\mu_\ell}{p_1\cdots p_{n+1}}h(0)<0$,
     which contradicts the fact that $\displaystyle{\frac{\partial c_2}{\partial \lambda_1}(0, 1)=1>0}$ by \eqref{Eq::Diff of c2}.

    By Corollary \ref{cor::Sign of Eigenvalues}, $M(\gamma)$ has at least 3 positive eigenvalues and 2 negative eigenvalues. Also, by Lemma \ref{lem::Sign Changes} and the fact that $\displaystyle{\lim\limits_{\mu\to 0^-}\frac{\partial c_1}{\partial \lambda_1}(\mu, 1)=+\infty}$, the Hessian $M(\beta)\mu+M(\gamma)$ has at least 2 positive eigenvalues and 2 negative eigenvalues when $\mu \in (-\infty,0)$. It follows that the real zero set of $f$ has at most 2 pieces, by Lemma \ref{lem::Cases when the number changes} and Lemma \ref{Lem::Number in Outer Chambers}.

     \item \textbf{$\Gamma_{\varepsilon}(\mathcal{A},B)$ has exactly 4 cusps.} By Lemma \ref{Lem::Structure of 4 cusps}, the depth of inner chambers in this case is at most two. We aim to show that the number of pieces in a depth-2 inner chamber matches that of an adjacent depth-1 chamber. Specifically, the boundary of a depth-2 inner chamber has at least 3 edges, each originating from a distinct curve segment. By Corollary \ref{cor::Sign of Eigenvalues} and Lemma \ref{lem::Sign Changes}, there are only two curve segments on which $M(\beta)\mu+M(\gamma)$ may have exactly 1 positive or 1 negative eigenvalue, across which the number of pieces could change. Thus, among the 3 edges on the boundary of a depth-2 inner chamber, there is an edge along which $M(\beta)\mu + M(\gamma)$ has at least 2 positive and 2 negative eigenvalues. By Lemma \ref{lem::Cases when the number changes}, this depth-2 inner chamber has the same number of pieces as the adjacent chamber across this edge, which, by Lemma \ref{Lem::Structure of 4 cusps}, has depth one. Hence, the depth-2 inner chamber has the same number of pieces as a depth-1 chamber. Consequently, the real zero set of $f$ contains at most 3 pieces.

      \item \textbf{$\Gamma_{\varepsilon}(\mathcal{A},B)$ has exactly 2 or 3 cusps.}  By Corollary \ref{Lem::Structure of 2 or 3 cusps} and \ref{Cor::Rough bound}, the real zero set of $f$ has at most 3 pieces.

      \item \textbf{$\Gamma_{\varepsilon}(\mathcal{A},B)$ has exactly 0 or 1 cusps.}  By Lemma \ref{Lem::One Cusp} and \ref{Lem::Number in Outer Chambers}, the real zero set of $f$ has at most 2 pieces. \qed
\end{itemize}

\end{document}